\documentclass[11pt]{article}
\catcode`\@=11 \catcode`\@=12
\usepackage[utf8]{inputenc}
\usepackage[T1]{fontenc}
\usepackage{amsmath}
\usepackage{amssymb}
\usepackage{amsthm}
\usepackage{oldgerm}
\usepackage{multicol}
\usepackage{verbatim}
\usepackage{color}
\usepackage{ulem}
\usepackage{enumitem}
\usepackage{mathrsfs}
\usepackage{graphicx}

\usepackage[pdftex,colorlinks,urlcolor=blue,pdfstartview=FitH,
]{hyperref}

\makeatletter
\newcommand{\mylabel}[2]{#2\def\@currentlabel{#2}\label{#1}}
\makeatother

\newcommand{\Rm}{\mathbb{R}}

\newcommand{\mS}{\ensuremath{\mathcal{S}}}

\newcommand{\vs}{\vspace{.2cm}}

\newtheorem{lem}{Lemma}[section]

\newtheorem{thm}{Theorem}

\newtheorem{prop}[lem]{Proposition}

\def\proof {\noindent{\sc{Proof. }}}
\def\qed {\mbox{}\hfill {\small \fbox{}} \\}
\def\lto{\longrightarrow}
\def\lmto{\longmapsto}

\date{}

\begin{document}

\begin{center}
	\begin{huge}
		{\bf Normal form near orbit segments of convex Hamiltonian systems.}\\
	\end{huge}
	\vs
-----
\vs
\end{center}
\begin{small}
\begin{multicols}{2}
\noindent
Shahriar Aslani
\footnote{
	\includegraphics[height=0.3cm]{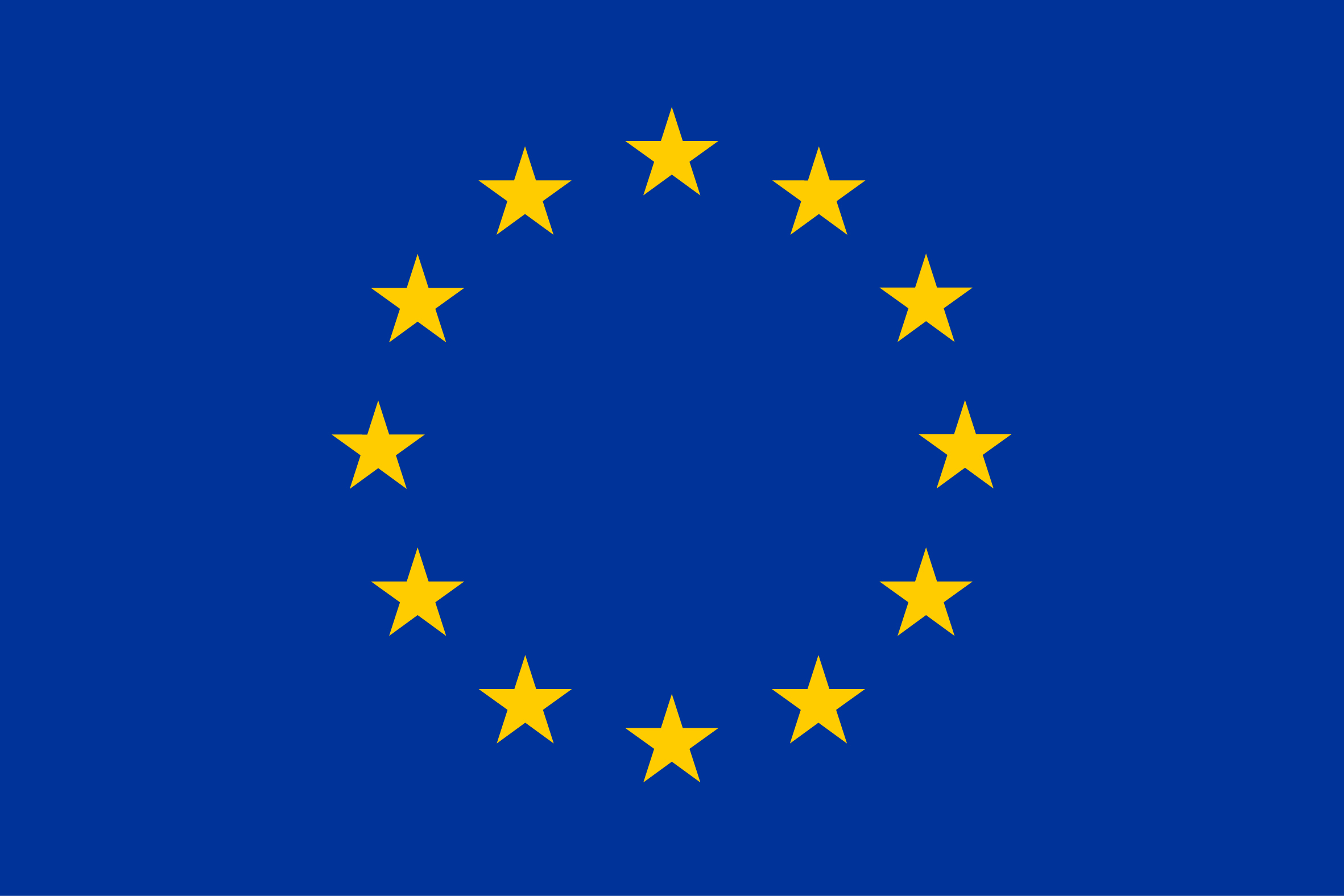} 
	This author has received funding from the European Union’s Horizon 2020 research and innovation programme under the Marie Skłodowska-Curie grant agreement No 754362}
\\
PSL Research University,\\
\'Ecole Normale Sup\'erieure,\\
DMA (UMR CNRS 8553)\\
45, rue d'Ulm\\
75230 Paris Cedex 05,
France\\
\texttt{shahriar.aslani@ens.fr}\\

\noindent
Patrick Bernard
\footnote{Universit\'e  Paris-Dauphine}\\
PSL Research University,\\
\'Ecole Normale Sup\'erieure,\\
DMA (UMR CNRS 8553)\\
45, rue d'Ulm\\
75230 Paris Cedex 05,
France\\
\texttt{patrick.bernard@ens.fr}\\

\end{multicols}
\end{small}
\vs

%
%
%
%
%
%
%
%
%
%
%
%
%
%
%

In the study of  Hamiltonian systems on cotangent bundles, it is natural to perturb Hamiltonians  by adding
potentials (functions depending only on the base point). This led to the definition of Ma\~né genericity: A property is generic if, given a  Hamiltonian $H$,
 the set of potentials $u$ such that $H+u$ satisfies the property is generic. This notion is mostly used in the context of Hamiltonians which are convex in $p$, in the sense that $\partial ^2_{pp} H$ is positive definite at
 each point. We will also restrict our study to this situation.

There is a close relation between perturbations of Hamiltonians by a small additive potential
and perturbations by a positive factor close to one. Indeed, the Hamiltonians
$H+u$ and $H/(1-u)$ have the same level one energy surface, hence  their dynamics on this energy surface are
reparametrisation of each other, this is the Maupertuis principle. This remark is particularly relevant when $H$ is homogeneous in the fibers (which corresponds to Finsler metrics) or even fiberwise quadratic
(which corresponds to Riemannian metrics). In these cases, perturbations by potentials of the Hamiltonian
correspond, up to parametrisation, to conformal perturbations of the metric.

One of the widely studied aspects is to understand to what extent the return map associated to a periodic orbit can be modified by a small perturbation.
This kind of question depend strongly on the context in which they are posed. Some of the
most studied contexts are, in increasing order of difficulty, perturbations of general vector fields, perturbations of Hamiltonian systems inside the class of Hamiltonian systems, perturbations of Riemannian metrics 
inside the class of Riemannian metrics, Ma\~né perturbations of convex Hamiltonians.
It is for example well-known that each vector field can be perturbed to a  vector field with only hyperbolic periodic orbits, this is part of the Kupka-Smale theorem, see \cite{K63,S63}. There is no such result in the context of Hamiltonian vector fields, but it remains true that each Hamiltonian can be perturbed to a Hamiltonian
with only non-degenerate periodic orbits (including the iterated ones), see \cite{CR1,CR2}.
The same result is true in the context of Riemannian metrics: every Riemannian metric can be perturbed to 
a Riemannian metric with only non-degenerate closed geodesics, this is the bumpy metric theorem, see \cite{KT,A82,A68}.
The question was investigated only much more recently in the context of Ma\~ né perturbations of convex Hamiltonians, see \cite{O08,RR}.
It is proved in \cite{RR} that the same result holds : If $H$ is a convex Hamiltonian  and $a$ is a regular value of $H$, then there exist arbitrarily small potentials $u$ such that all  periodic orbits
(including iterated ones)  of $H+u$ at energy $a$ are non-degenerate. The proof given in \cite{RR}
is actually rather similar to the ones given in papers on the perturbations of Riemannian metrics.
In all these proofs, it is very useful to work in appropriate coordinates around an orbit segment.
In the Riemannian case, one can use the so-called Fermi coordinates.
In the Hamiltonian case, appropriate  coordinates are considered in  \cite[Lemma 3.1]{RR} itself taken from \cite[Lemma C.1]{FR2}. However, as we shall detail below, the proof of this Lemma in \cite{FR2}, Appendix C, is incomplete, and the statement itself is actually wrong.
Our goal in the present paper is to state and prove  a corrected version of this normal form Lemma. Our proof is different from the one outlined  in \cite{FR2}, Appendix C.  In particular, it
is purely Hamiltonian and does not rest on the results of \cite{LN} on Finsler metrics, as \cite{FR2} did. Although our normal form
is weaker than the one claimed in \cite{RR}, it is actually sufficient to prove the main results of \cite{RR,LRR}, as we shall explain after the statement of Theorem 1,
and probably also of the other works using \cite[Lemma C.1]{FR2}.

\section{Introduction.}
When studying Ma\~né generic properties of convex Hamiltonians, a natural group of changes of coordinates to consider is the group of symplectic diffeomorphisms preserving the vertical fibration, \textit{i.e.} those symplectic diffeomorphisms of the form $\psi(q,p)=(\varphi(q), G(q,p))$.  We call such transformations \textit{fibered}.
It is well known that such a diffeomorphism is symplectic if and only  if the second coordinate $G(q,p)$ is of the form
$$
G(q,p)=\alpha_q+ p\circ (d\varphi_q)^{-1}
$$
for some closed one-form $\alpha$. We will say that the fibered symplectic diffeomophism $\psi$ is \textit{homogeneous} if 
 it preserves the zero section, which implies that there exists a diffeomorphism $\varphi$ of the base such that 
$$
\psi(q,p)=(\varphi(q), p\circ (d\varphi_q)^{-1}).
$$
We will say that the above diffeomorphism is  \textit{vertical} if it is of the form $(q,p)\lmto (q,p+\alpha_q)$, that is if it preserves the first coordinate.
Each fibered symplectic diffeomorphism is thus the composition of a vertical and of a homogeneous symplectic diffeomorphisms (in any order). 
Note that fibered symplectic diffeomorphisms preserve convexity, since their restrictions to fibers is affine.

If $\psi$ is a fibered symplectic diffeomorphism and $g(q)$ is a potentiel, then
 $(H\circ \psi)+g=(H+g\circ \varphi^{-1})\circ \psi$. So if a property is invariant under symplectic diffeomorphisms, for example having only non-degenerate periodic orbits on a given energy level,
 then this property is satisfied by $H+g$ for arbitrarily small $g$ if and only if it is satisfied by
 $g+H\circ \psi$ for arbitrarily small $g$.
 
 These considerations allow us to enlarge the group of transformations considered in \cite[Lemma C.1]{FR2}, where 
 only homogeneous symplectic  diffeomorphisms are considered.
 Although it is natural to restrict to homogeneous transformation in the case where $H$ is homogeneous, there is no reason to do so in general, and we will actually see that it is not possible to do so : Allowing vertical symplectic transformations is necessary to obtain a nice normal form.

Since all considerations are local, we will always work on the manifold $\Rm^{1+d}$ and on its cotangent bundle
$\Rm^{1+d}\times (\Rm^{1+d})^*$. We will use the notations $q=(q_0,q_*)\in \Rm\times \Rm^d$ and similarly
$p=(p_0,p_*)\in \Rm^*\times \Rm^{d*}$. We denote by $e_0, e_1, \ldots, e_d$ both the standard base of $\Rm^{1+d}$
and $\Rm^{(1+d)*}$.

If $(Q,P)$ is an orbit segment such that $\dot Q(0)\neq 0$, then there exists a local diffeomorphism
$\varphi$ of the base which sends the orbit segment $Q(t), t\in[-\delta, \delta]$ to the straight line segment 
$te_0, t\in[-\delta, \delta]$.
Once this reduction has been performed, we  only consider 
fiber-preserving symplectic diffeomorphims $\psi$ which have the property that their horizontal component $\varphi$  is the identity on the 
segment $[-\delta, \delta]e_0$.
We call such diffeomorphisms \textit{admissible}.  This implies that 
\begin{equation}\label{eq-e0}
d\varphi_{te_0} \cdot e_0=e_0
\quad, \quad
(p \circ (d\varphi_{te_0})^{-1})_0=p_0,
\end{equation}
where  the second equality is obtained by applying the linear form $p\circ d\varphi^{-1}_{te_0}$ to the first equality,
and where $p_0$ is the first coordinate $p_0=p\cdot e_0$. In other words, the first component 
$p_0$ of the momentum is not changed by applying an admissible homogeneous symplectic diffeomorphism.

\begin{thm}\label{thm-1}
Let $\underline H: T^*\Rm^{1+d} \lto \Rm $ be a smooth Hamiltonian convex in $p$ in the sense that  $\partial^2_{pp}H$ 
is positive definite at each point. Let  $(\underline Q(t), \underline P(t))$ be an orbit of $\underline H$ such that 
$\dot {\underline Q}(0) \neq 0$.
Then there exists a smooth  local fibered symplectic diffeomophism $\psi$ and $\delta >0$  such that the new Hamiltonian 
$H=\underline H\circ \psi$ and the new orbit $(Q(t), P(t))=\psi^{-1}(\underline Q(t), \underline P(t))$ satisfy, for all 
$t\in [0, \delta]$, 
\begin{eqnarray}
\label{eq-Q}
Q(t) &=& te_0,\\
\label{eq-P}
P(t) &=& 0,\\
\label{eq-p0p*}
\partial^2 _{p_0p_*}H(te_0,P(t)) &=&0,\\ 
\label{eq-qp*}
\partial ^2_{qp_*} H (te_0, P(t)) &=& 0,\\
\label{eq-p*p*}
\partial^2 _{p_*p_*}H(te_0,P(t)) &=&Id.
\end{eqnarray}
If  $\partial^2_{q_*p_0} H (te_0, 0)$  does not identically vanish, then this quantity  can not be reduced to zero by applying an admissible change of coordinates preserving the other equalities.
\end{thm}

Note that equality (\ref{eq-P}) can obviously not be obtained using only homogeneous changes of coordinates,
since they preserve the zero section.
Moreover, it follows from  (\ref{eq-e0}) that the first coordinate $P_0(t)$ is invariant under the action 
of homogeneous admissible diffeomorphisms. 
As a consequence, it is not true that 
orbits of general convex  Hamiltonians can be reduced to  $(Q(t),P(t))=(te_0, e_0)$ 
by such diffeomorphisms, 
as is claimed in 
\cite[Lemma C.1]{FR2}. This becomes possible (and easy)
once vertical changes of coordinates are allowed,   as we shall verify below.

Contrarily  to \cite[Lemma C.1]{FR2},  the last claim of the theorem  implies that (\ref{eq-qp*}) can't be strengthened to
$\partial^2_{qp}H (te_0, 0)=0$, 
even if   vertical symplectic diffeomorphisms, in addition to homogeneous ones, are permitted.
See however Section \ref{sec-h} where this equality is proved in the homogeneous case.

If  $H$ satisfies all the conclusions of Theorem \ref{thm-1}, it is of the form 
\begin{equation}\label{eq-h}
H(q,p)=f(q)+w(q)p_0+\frac{1}{2}  a(q_0) p_0^2 
+\frac{1}{2} \langle  p_*, p_*\rangle + O_3(q_*, p),
\end{equation}
where $f$ and $w$ are smooth functions from $\Rm^{1+d}$ to $\Rm$ satisfying  $f(q_0,0)=f(0,0)$ and  $w(q_0,0)=1$ for each $q_0\in [0, \delta]$. Setting $\tilde H(q,p):= (H(q,p)-H(0))/w(q)$, we thus have 
\begin{equation}\label{eq-th}
\tilde H(q,p)= \tilde f(q)+p_0+\frac{1}{2}  \tilde a(q_0) p_0^2 
+\frac{1}{2} \langle  p_*, p_*\rangle + O_3(q_*, p),
\end{equation}
with $\tilde f(q_0,0)=0$.
This means that the additional conclusion $\partial^2_{qp}H(te_0,0)=0$ can be achieved provided 
we translate $H$ so that our orbit has energy $0$ (which does not change anything to the dynamics),
and then multiply $H$ by a function of $q$ (which is a reparametrization of the dynamics on the
energy surface $\tilde H^{-1}(0)=H^{-1}(H(0))$). 
This fact seems sufficient to derive   most of the applications of \cite[Lemma C.1]{FR2} existing in the literature. 
 
Let us illustrate for example how the results of \cite{RR} can be obtained.
We denote by $E_t$ the space $\{q_0=t\}\cap \{p_0=0\}$, which projects isomorphically to 
$\Rm^{d}\times \Rm^{d*}$.
Let $L:E_0\lto E_{\delta}$ be the differential at zero of the transition maps between  the sections 
$\{q_0=0\}\cap \{\tilde H=0\}$ and $\{q_0=\delta\}\cap \{\tilde H=0\}$, seen as a  symplectic $2d\times 2d$ matrice. Since the dynamics of $H$ and $\tilde H$
in restriction to $\tilde H^{-1}(0)$ are reparametrizations of each other, they have the same transition map
hence $L_H= L_{\tilde H}$. 
The main statement of \cite{RR} is that, if $X$ is a dense set of symplectic $2d\times 2d$ matrices,  there exist 
arbitrarily  small potentials $g(q)$ such that $L_{H+g}\in X$.
This statement is proves in \cite{RR} for Hamiltonian having the form (\ref{eq-th}) above, and the normal form is invoked,
to reduce each Hamiltonian to this form. However, as we have explained, one can only obtain the normal form (\ref{eq-h}) in general.
The missing step is to deduce the statement for $H$ from the statement for $\tilde H$, which turns out to be easy :
Applying the statement to $\tilde H$ gives  small potentials $\tilde g(q)$ such that 
$L_{\tilde H+\tilde g}\in X$. We now observe that $\tilde H +\tilde g=(H-H(0)+w\tilde g)/w$, hence 
$$
L_{H+w\tilde g }=L_{\tilde H+\tilde g}\in X.
$$
Since the function $w$ depends only on $q$, $w\tilde g$ is a potential, which can be made arbitrarily small by taking $\tilde g$ arbitrarily small. We have proved the statement for $H$.

\section{Proof of the Normal Form.}

We will always work in coordinates such that (\ref{eq-Q}) holds, and consider only admissible changes of coordinates (\textit{i.e.} changes of coordinates whose horizontal component fixes the axis $\Rm e_0$).
Our proof is purely Hamiltonian, and does not rest on \cite[Lemma 3.1]{LN}. Actually, a small modification of the proof also implies this Lemma, as will be explained in the next section.
We will apply several  admissible diffeomorphisms. At each step,
we will denote by $\underline H$ the initial Hamiltonian  and by $H=\underline H \circ \psi$ the transformed Hamiltonian.
In the matrix computations below, we most of the time consider the momenta $p$ as line matrices hence denote by $pM$
what might also be denoted by $M^tp$.

\noindent
\textsc{proof of (\ref{eq-P}).}
Let $\underline P_0(t)$ be the first component of $\underline P(t)$ (the orbit before the change of coordinates).
We consider a function $v(t) : \Rm\lto \Rm$ such that  $v'=\underline P_0$
and the function 
$u(q_0, q_*):=v(q_0)+\underline P_*(q_0)\cdot q_*.
$
We have $du _{te_0}=\underline P(t)$, hence applying the vertical diffeomorphism $\psi (q,p)=(q, p+du_q)$,
the new orbit $(Q(t), P(t))=\psi^{-1}(\underline Q(t), \underline P(t))$ satisfies $P(t)=0$. \qed

\noindent
\textsc{proof of (\ref{eq-p0p*}).}
We assume that (\ref{eq-P}) and (\ref{eq-Q}) are already satisfied for $\underline H$, and prove that 
 (\ref{eq-p0p*}) can be obtained by a further change of coordinates.
We consider a base diffeomorphism of the form
$\varphi (q_0,q_*)= (q_0+l(q_0)\cdot q_*, q_* )$,
where $q_0\lmto l(q_0)$ is a smooth map with values in $\Rm^{d*}$.
The corresponding homogeneous diffeomorphism satisfies 
\begin{equation}
\label{eq-psi}
\psi: (q_0,0,p_0,p_*)\lmto (q_0, 0, p_0, p_*+p_0 l(q_0)).
\end{equation}
We then have 
$\partial _{p_*}(\underline H\circ \psi)_{(q_0e_0, 0)}=
\partial _{p_*}\underline H _{(q_0e_0, 0)}
$
hence
$$
\partial^2_{p_0 p_*}H_{(q_0e_0,0)}=\partial^2_{p_0 p_*}(\underline H\circ \psi)_{(q_0e_0,0)}=\partial^2 _{p_0p_*}\underline H _{(q_0e_0,0)}+
\partial^2_{p_*p_*}\underline H_{(q_0e_0,0)}\cdot l(q_0).
$$ 
We obtain (\ref{eq-p0p*}) by choosing 
$$
l(q_0):= -(\partial^2_{p_*p_*}\underline H_{(q_0e_0,0)})^{-1} \cdot
\partial^2 _{p_0p_*}\underline H _{(q_0e_0,0)}.
$$
Observe that  $\partial^2_{p_*p_*}H_{(q_0e_0,0)}$ is invertible because
$\partial^2_{pp}H$ is positive definite at each point.
 \qed

\noindent
\textsc{proof of (\ref{eq-qp*}).}
This equality can be obtained by a further homogeneous change of coordinates preserving (\ref{eq-Q}) and (\ref{eq-p0p*}). We could assume (\ref{eq-P}), but, keeping in mind another application in Section \ref{sec-h}, we  only  make the slightly  more general assumption that 
$\underline P(t)\equiv \underline P(0)=(P_0,0)$ for some constant $P_0$.
We consider the vector field
$$
\underline V(q):= \partial _p \underline H(q, P(0))
$$
on $\Rm^{d+1}$.
We will  apply a variant of the Flow Box Theorem to the vector field $\underline V$.
More precisely, we consider the diffeomorphism $\varphi(q)= (q_0, \phi(q_0,q_*))$, defined in a neighborhood of 
$[-\delta , \delta]\times \{0\}$ in such a way that $q_*\lmto  \phi(t, q_*)$ is the transition map along the orbits of $\underline V$
between the sections $\{q_0=0\}$ and $\{q_0=t\}$. In other words, $\varphi(q_0, q_*)=\Phi^{q_0}(0,q_*)$,
where $\Phi^t$ is the flow of 
the reparametrized vector field $\underline V(q)/\underline V_0(q)$ ($\underline V_0$ is the first coordinate of $\underline V$).
 It is a smooth diffeomorphism near $[-\delta, \delta]\times \{0\}$,
and $\underline V=\underline V_0 \;\varphi_{\sharp} e_0$, where $\varphi_{\sharp}e_0$ is the forward image of the constant  vector field $e_0$.

Since $\varphi$ is fixing the axis $\Rm e_0$ and preserving $q_0$,
the associated  homogeneous diffeomorphism $\psi$  preserves (\ref{eq-Q}) and (\ref{eq-p0p*}).
Moreover, $\psi (q, P(0))=(\varphi(q), P(0)).$
Denoting as usual $H:= \underline H\circ \psi$, and $V(q)=\partial_p H(q,P(0))$, we have  
$
\varphi_{\sharp} V=\underline V =\underline V_0 \;\varphi_{\sharp} e_0,
$
hence $V=(\underline V_0\circ \varphi^{-1})e_0$, and  $V_*=0$, or in other words $\partial_{p_*} H(q,P(0))=0$
for all $q$. Differentiating with respect to $q$ gives (\ref{eq-qp*}). 
\qed

\noindent
\textsc{proof of (\ref{eq-p*p*}).}
We assume that the equations (\ref{eq-P}) to (\ref{eq-qp*}) initially hold.
We will obtain (\ref{eq-p*p*}) by an admissible (usually not homogeneous) transformation preserving all these equalities. This transformation will be decomposed  into first a homogeneous tranformation and second a vertical
transformation none of which  preserve (\ref{eq-qp*}).

The first step consists of applying the homogeneous change of coordinates $\psi$ associated to a diffeomorphism
of the form
$$
\varphi (q_0,q_*)=(q_0, M(q_0)\cdot q_*),
$$
where $M(t)$ is a $d\times d$  invertible matrix depending smoothly on $t$.
The matrix of the differential of $\varphi$ is 
$$
D(q)=
\begin{bmatrix}
1 & 0\\ M'(q_0)q_* & M(q_0)
\end{bmatrix},
\quad 
D^{-1}(q)=
\begin{bmatrix}
1 & 0\\- M^{-1}(q_0)M'(q_0) q_* & M^{-1}(q_0)
\end{bmatrix},
$$
where $M'(q_0)$ is the derivative.
We thus have
$$
\psi(q,p)=\big(q_0, M(q_0) q_*, p_0-p_* M^{-1}(q_0)M'(q_0)  q_*
,p_* M^{-1}(q_0)\big).
$$
The Hamiltonian in original coordinates is of the form
$$
\underline H (q,p)=
\underline H(q,0)+ \underline v(q) p_0+ \frac{1}{2} \underline a(q_0) p_0^2 + \frac{1}{2} \langle  p_*\underline A(q_0), p_*\rangle + O_3(q_*, p),
$$
with $\underline v(q)=\partial_{p_0}\underline H(q,0)$, $\underline a(q_0)=\partial^2_{p_0p_0} \underline H(q_0e_0, 0)$,
$\underline A(q_0)=\partial^2_{p_*p_*} \underline H(q_0e_0,0)$.
We compte
\begin{align*}
H(q,p)=\underline H\circ \psi (q,p) &=
H(q,0)+v(q)p_0-p_*M^{-1}(q_0)M'(q_0) q_*+\frac{1}{2}  a(q_0) p_0^2 \\
& +\frac{1}{2} \langle  p_*M^{-1}(q_0)\underline A(q_0), p_*M^{-1}(q_0)\rangle + O_3(q_*, p).
\end{align*}
We get (\ref{eq-p*p*}) provided $M(q_0)M^t(q_0)=\underline A(q_0)$ for each $q_0$. We could for example take 
$M(q_0)=\underline A^{1/2}(q_0)$ for each $q_0$ (remember that $\underline A(q_0)$ is positive definite for each $q_0$).
However, the unavoidable  apparition of the term $p_*M^{-1}(q_0)M'(q_0) q_*$ means that (\ref{eq-qp*}) have been destroyed. In order to be able to restore it by a vertical change of coordinates, we need a better choice for  $M(q_0)$ :

\begin{lem}
	We can choose $M(q_0)$ in such a way that 
	$$
	H(q,p)=
	H(q,0)+v(q)p_0-p_*B(q_0) q_*+\frac{1}{2}  a(q_0) p_0^2 
	+\frac{1}{2} \langle  p_*, p_*\rangle + O_3(q_*, p),
	$$ 
	where $B(q_0)$ is symmetric for each $q_0$.
\end{lem} 

\proof
We need the matrix $M(t)$ to satisfy the two conditions that $M(t) M^t(t)=\underline A(t)$ and $B(t):=M'(t)M^{-1}(t)$ is 
symmetric. Derivating the first condition, we get $M'M^t +M(M')^t =\underline A'$.
Expressing this relation in terms of $B$, we obtain the equation
$
BMM^t+MM^t B=\underline A'
$,
 which can be rewritten 
 $$B\underline A+\underline A B=\underline A'.
 $$
 We claim that the linear map $L_P:X\lmto PX+XP$  from $\mS_d(\Rm)$ (the set of symmetric matrices) to itself is an isomorphism if $P$
is a positive definite symmetric matrix.
Indeed, if $P$ is diagonal, then $(PX+XP)_{ij}=(P_{ii}+P_{jj})X_{ij}$, with $P_{ii}+P_{jj}>0$. 
In general, $P=ODO^{t}$ for some orthogonal matrix $O$, 
and then $L_P(X)=OL_D(O^{t}XO)O^{t}$, hence $L_P$ is also an isomorphism.
The map $P\lmto L_P$ is smooth, hence so is the map $P\lmto (L_P)^{-1}=:L^{-1}_P$.
We deduce that the equation $R(t)\underline A(t)+\underline A(t)R(t)=\underline A'(t)$ has a unique symmetric solution 
$$
B(t)=L^{-1}_{\underline A(t)}(\underline A'(t)),
$$
this solution depends smoothly on $t$.
With this $B(t)$, we consider the solution $M(t)$ of the differential equation 
$$
M'(t)=B(t)M(t)
$$
with initial condition $M(0)=\underline A^{1/2}(0)$. 
By definition of $B$, we have
$$
\underline A'(t)=M'(t)M^{-1}(t)\underline A(t)+\underline A(t)M'(t)M^{-1}(t).
$$
Seeing this as a differential equation in $A$, we just need to check that  $M(t)M^t(t)$ solves this equation
to deduce that $A(t)=M(t)M^t(t)$.
This follows from the simple computation
$$
(MM^t)'=M'M^t+M(M^t)'=M'M^{-1}(MM^t)+(MM^t)M'M^{-1},
$$
where we have used that $(M^t)'=M^tM'M^{-1}$, which holds because $M'M^{-1}$ is symmetric.
\qed

The second step consists of applying the vertical change of coordinates 
$$
\Psi :(q,p)\lmto (q,p+du_q)
$$
with $u(q)=\langle B(q_0)q_*, q_*\rangle/2$, so that $du_q=(*,B(q_0)q_*)$. It is then a direct computation that  
$$
H\circ \Psi (q,p)=f(q)+w(q)p_0+\frac{1}{2}  a(q_0) p_0^2 
+\frac{1}{2} \langle  p_*, p_*\rangle + O_3(q_*, p),
$$
for some smooth functions $f$ such that $f(q_0,0)=H(0,0)$, and $w$ such that $w(q_0,0)=1$.\qed

We now prove the last statement of the theorem, about the impossibility of achieving the additional condition
 $\partial^2_{q_*p_0}H=0$.
We shall only consider admissible diffeomorphisms which preserve (\ref{eq-Q}), (\ref{eq-P}) and  (\ref{eq-p0p*}). 
Every fibration preserving symplectic diffeomorphism preserving (\ref{eq-Q}) and  (\ref{eq-P}) is the composition
of a homogeneous and a vertical diffeomorphisms each of which preserve (\ref{eq-Q}) and (\ref{eq-P}).

Let us first observe that $\partial^2_{q_*p_0}H(te_0,0)$ can't be changed by applying a vertical diffeomorphism
preserving (\ref{eq-P}). Such a diffeomorphism is of the form $\psi(q,p)=(q,p+du_q)$ for some smooth function $u$
satisfying $du_{te_0}=0$, hence in particular $u$ is constant on $\Rm e_0$.
Then, 
$$
\partial_{p_0}(H\circ \psi) (q,p)=\partial_{p_0}H (q,p+\partial_{q}u(q))
$$
and
\begin{align*}
\partial^2_{q_*p_0}(H\circ \psi)(te_0,0)
&=\partial^2_{q_*p_0}H((te_0,0)+\sum_{i=0}^ d\partial^2_{p_ip_0}H(te_0,0)\partial^2_{q_*q_i}u(te_0)\\
&= \partial^2_{q_*p_0}H((te_0,0)+\partial^2_{p_0p_0}H(te_0,0)\partial^2_{q_*q_0}u(te_0)
=\partial^2_{q_*p_0}H(te_0,0).
\end{align*}
In this computation, we have used first that 
$\partial^2_{p_*p_0}H(te_0,0)=0$, and then that 
$\partial^2_ {q_*q_0}u(te_0)=\partial^2_ {q_0q_*}u(te_0)=0$, which holds since   $\partial_{q_*}u(te_0)$ is identically zero.

We now consider the action of homogeneous admissible transformations.

\begin{lem}
	The homogeneous  symplectic diffeomorphism $\psi$ associated to $\varphi$ preserves (\ref{eq-Q}) and (\ref{eq-p0p*}) if and only if the matrix of the differential  of $\varphi$  has the following  $1+d$ block form for each $t\in [0,\delta]$:
	$$
	D\varphi (te_0,P(t))= 
	\begin{bmatrix} 1 & 0\\ 0& *\end {bmatrix}.
	$$
\end{lem}

\proof
Since $\varphi$ is admissible, the matrix $D\varphi$ along the orbit  has the triangular block  form
$$D(t):=D\varphi (te_0,P(t))= 
\begin{bmatrix} 1 & b(t)\\ 0& B(t)\end{bmatrix},
\quad 
D^{-1}(t)= 
\begin{bmatrix} 1 & -b(t)B^{-1}(t)\\ 0& B^{-1}(t)\end{bmatrix},
$$
hence 
$$
\psi(q,p)=\big(\varphi(q), p_0, p_*B^{-1}(t)-p_0b(t)B^{-1}(t))\big).
$$
We have 
$$
\partial^2_{pp} (H\circ \psi)(te_0, P(t))
=D^{-1}(t) \partial^2_{pp}H_{\psi (te_0,P(t))}D^{-1t}(t).
$$
In matrix form,
\begin{align*}
\partial ^2_{pp}(H\circ \psi)
&=
\begin{bmatrix} 1 & -b(t)B^{-1}(t)\\ 0& B^{-1}(t)\end {bmatrix}
\begin{bmatrix} \partial^2_{p_0p_0}H & 0\\ 0& \partial^2_{p_*p_*}\end {bmatrix}
\begin{bmatrix} 1 & 0\\ -B^{-1t}(t) b^t(t)& B^{-1t}(t)\end {bmatrix}\\
&=\begin{bmatrix} * & -b(t)B^{-1}(t)A(t)B^{-1t}\\ B^{-1}(t) A(t) B^{-1t}(t)b^t& B^{-1}(t)A(t) B^{-1t}\end {bmatrix}.
\end{align*}
This matrix is block diagonal (which is equivalent to $\partial^2_{p_0p_*}(H\circ \psi)=0$)
if and only if $b(t)=0$.
\qed

Let now $\underline H$ be a Hamiltonian satisfying (\ref{eq-Q}) to (\ref{eq-p0p*}), 
and $\psi$ be the homogeneous transformation associated to the diffeomorphism  $\varphi$.
We assume that $H=\underline H\circ\psi$ still satisfies (\ref{eq-Q}) to (\ref{eq-p0p*}).
We denote by $a(q)$ the first coordinate of $\varphi$.
In view of the previous Lemma, we have the block diagonal form for $t\in [0, \delta]$,
$$D\varphi_{(te_0,0)}=\begin{bmatrix} 1&0\\0& Z(t)\end{bmatrix},
$$
in particular $\partial_q a(q_0,0)=e_0$, 
$\partial_{q_*} a(q_0,0)=0$ for each $q_0 \in[0, \delta]$ thus
$\partial^2_{q_0q_*} a (q_0,0)=0$. 
It is convenient for the following computations to denote $V(q):= \partial_q H(q,0)$,
$\underline V(q):= \partial_q \underline H(q,0)$,
so that 
$$
\underline V(\varphi(q))=d\varphi_q \cdot V(q)
$$
and, focusing on the first coordinates, $\underline V_0(\varphi(q))=\partial_q a(q)\cdot V(q).$
Differentiating with respect to $q_*$ at the point $te_0$, $t\in [0, \delta]$ yields :
\begin{align*}
\partial_{q_*}\underline V_0(te_0)\cdot Z(t)
&=
\partial^2_{q_*q}a(te_0) \cdot V(te_0)
+ \partial_q a(te_0)\cdot \partial_{q_*}V(te_0)\\
&=\partial^2_{q_*q}a (te_0)\cdot e_0 + \partial_{q_*}V_0(te_0)
= \partial^2_{q_*q_0}a(te_0)+\partial_{q_*}V_0(te_0)\\
&=\partial_{q_*}V_0(te_0)
\end{align*}
Since $Z(t)$ is invertible for each $t\in[0,\delta]$,  
we have proved that $\partial^2_{q_*p_0}H(te_0,0)=\partial_{q_*}V_0(te_0)=0$ if and only if
$\partial^2_{q_*p_0}\underline H(te_0,0)=\partial_{q_*}\underline V_0(te_0)=0$.
If $\underline H$ did not satisfy this condition from the start, then neither does $H$.
\qed

\section{The Homogeneous case}\label{sec-h}

We explain here for completeness how the arguments given above also imply the following statement, which is equivalent  to 
\cite[Lemma 3.1]{LN}.

\begin{prop}\label{thm-2}
	Let $\underline H: T^*\Rm^{1+d}\lto \Rm $  be a  Hamiltonian positively homogeneous in the fibers, smooth and positive outside of the zero section, and such that $\partial_{pp} (H^2)$ is positive definite at each point outside of the zero section.
	Let  $(\underline Q(t), \underline P(t))$ be an orbit of $\underline H$ such that 
	$\dot {\underline Q}(0) \neq 0$.
	Then there exists a smooth  local fibered homogeneous  symplectic diffeomophism $\psi$ and $\delta >0$  such that the new Hamiltonian 
	$H=\underline H\circ \psi$ and the new orbit $(Q(t), P(t))=\psi^{-1}(\underline Q(t), \underline P(t))$ satisfy, the 
	equalities (\ref{eq-Q}), (\ref{eq-p0p*}), (\ref{eq-qp*})
	for each $t\in [0, \delta]$, as well as
	\begin{eqnarray}
	\label{eq-PH}
	P(t) &= & (P_0(0),0),\\
	\label{eq-qp}
	\partial ^2_{qp} H (te_0, P(t)) &=& 0.
	\end{eqnarray}
\end{prop}

\proof
We work in coordinates where (\ref{eq-Q}) hold and apply admissible homogeneous transformations.

The first component
$P_0(t)$ of the momentum is  independent of $t$, and non zero.
Indeed, we have, using the Euler relation 
$$
P_0(t)=P(t)\cdot e_0 =P(t)\cdot \partial _p H(te_0, P(t))=a H(te_0, P(t)),
$$
where $a$ is the degree of homogeneity.
Then,
$$P_0(t)=a H(te_0, P(t))=a H(0,P(0))=P_0(0),
$$
where we have used the preservation of $H$ along the orbit. 
Moreover, the constant $P_0(t)$ as well as the energy $H(te_0, P(t))$ are non-zero otherwise the orbit would be constant.

We now apply the same homogeneous diffeomorphism  as in the proof of (\ref{eq-p0p*}) above, with
$l(q_0)=\underline P_*(t)/P_0$, where $P_0=\underline P_0$ is the first component of the momentum.  In view of (\ref{eq-psi}), we get (\ref{eq-PH}).

Then  (\ref{eq-p0p*}) automatically holds : We have 
$$\partial_{p_*} H (te_0, (P_0,0))=\partial_{p_*}  H (te_0,P(t))=\dot Q_*(t)=0.
$$  
By homogeneity of the map $p\lmto \partial_{p_*} H (te_0, p)$, we deduce that 
$\partial_{p_*} H (te_0, (sP_0,0))=0$ for each $s>0$, hence that 
$\partial^2_{p_0p_*} H(te_0, (p_0,0))=0$.

The equality (\ref{eq-qp*}) can then be obtained by applying a homogeneous change of coordinates, the proof of the previous section can  directly be applied.

The equality  $\partial_ {qp_0}H(te_0, P(t))=0$, hence 
(\ref{eq-qp}), 
then automatically holds.
Indeed, since $P(t)=P(0)$ is constant, we have 
$\partial_q H(te_0,(P_0(t),0))=\partial_qH(te_0,P(t))=\dot P(t)=0$. Using the homogeneity of the map 
$p\lmto \partial_qH(te_0,p)$, we deduce  that $\partial_q H(te_0,(sP_0,0))=0$ for each $s>0$,
hence that 
$\partial_{p_0q}H(te_0, (P_0,0))=0$.
\qed

\end{document}